\documentclass{amsproc}
\newif\ifpdf
    \ifx\pdfoutput\undefined
    \pdffalse 
    \else
    \pdfoutput=1 
    \pdftrue
    \fi

    \ifpdf
    \usepackage[pdftex]{graphicx}
    \else
    \usepackage{graphicx}
    \fi
\usepackage{amsmath,amssymb,amsthm,epsfig}

\newcommand\Z{\mathbb Z}

\newtheorem{theorem}{Theorem}[section]
\newtheorem{proposition}[theorem]{Proposition}
\newtheorem{lemma}[theorem]{Lemma}

\theoremstyle{definition}
\newtheorem{definition}[theorem]{Definition}

\input{epsf.tex}

\begin{document}
\ifpdf
    \DeclareGraphicsExtensions{.pdf, .jpg, .tif}
    \else
    \DeclareGraphicsExtensions{.eps, .jpg}
    \fi

\title{Metric properties of the lamplighter group as an automata group}
\author{Sean Cleary}
\address{Department of Mathematics,
The City College of New York,
New York, NY 10031}
\email{cleary@sci.ccny.cuny.edu}
\thanks{The first author acknowledges support from PSC-CUNY grant \#64459-0033}
\author{Jennifer Taback}
\address{Department of Mathematics,
Bowdoin College,
Brunswick, ME 04011}
\email{jtaback@bowdoin.edu}
\thanks{The second author acknowledges support from
NSF grant DMS-0437481}


\begin{abstract}
We examine the geometry of the Cayley graph  of the lamplighter
group with respect to the generating set rising from its
interpretation as an automata group due to Grigorchuk and Zuk.  We
find some metric behavior with respect to this generating set analogous
to the metric behavior in the standard group theoretic generating
set. The similar metric behavior
includes expressions for geodesic paths and families of
`dead-end' elements, which are endpoints of terminating geodesic
rays. We also exhibit some different metric behavior between these
two generating sets related to the existence of `seesaw' elements.
 \end{abstract}

\maketitle

\section{Introduction}
\label{sec:intro}

There are several generating sets of interest for the lamplighter
group $L$, which is the wreath product of the group of order two
with the integers. The standard presentation of $L$ arises from
this wreath product structure, namely
$$
L = \langle a,t | a^2, [t^{i}a t^{-i},t^{j}at^{-j}],  i,j \in \Z
\rangle.$$ We will refer to the set of elements $\{a,t\}$ as the
wreath product generating set of $L$.

\begin{figure}\includegraphics[width=3.5in]{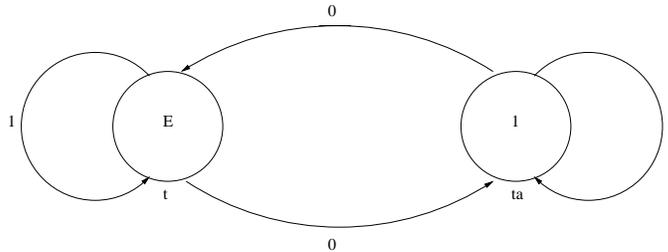}
\caption{The two state automaton realizing the generators of the
lamplighter group $t$ and $ta$.
\label{fig:automatonpic}}
\end{figure}

Automata groups form a rich class of groups with a number of
remarkable properties.  For example, Grigorchuk's groups of
intermediate growth are realized by automata groups, as described
by Grigorchuk, Nekrashevich and Sushchanski\u{\i}  in
\cite{automata}. Grigorchuck and {\.Z}uk show in \cite{zukspec}
that the lamplighter group $L$ can also be regarded as an
automata group. The automaton they use to realize  $L$ in
this way is shown in Figure \ref{fig:automatonpic}. The generating
set arising from this interpretation is $\{t, ta\}$, where $a$ and
$t$ are the wreath product generators of $L$. We will refer to
the generating set $\{t, ta\}$ as the automata generating set of
$L$.
  The generator $t$
corresponds to the initial automaton with start state labelled `t'
and the generator $ta$ corresponds to the initial automaton with
start state labelled `ta'. Grigorchuk and {\.Z}uk compute the
spectral radius of $L$ with respect to this generating set; they
find remarkably that it is a discrete measure.

We explored the geometry of the Cayley graph of $L$ with respect
to the wreath product generating set $\{a,t\}$ in \cite{deadlamp}.
Although the perspective which gives rise to the automata
generating set is quite different from the wreath product
perspective, we show below that the  metric properties of the
group are  similar  in several ways with respect to either
generating set. In particular, we find that  {\em dead-end}
elements of arbitrary depth occur with respect to both generating
sets. With respect to the wreath product generating set there are
`seesaw' elements \cite{deadlamp}, which do not occur with respect
to the automaton generating set.  However, we  find and describe
{\em seesaw-like} behavior with respect to the automaton
generating set.  It is not known the extent to which the existence
of dead-end elements of arbitrary depth depends upon the choice of
generating set. We show here that the existence of seesaw words
does depend upon generating set.

\section{Normal forms and geodesics in $L$ }
\label{sec:lengths}

\begin{figure}\includegraphics[width=3in]{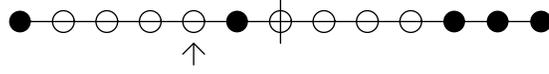}\\
\caption{The element $w=a_4 a_5 a_6 a_{-1} a_{-6} t^{-2}$ of $L$,
expressed as a configuration of illuminated bulbs and a cursor in
the position $-2$.  Solid circles represent bulbs which are
illuminated, open circles  represent bulbs which are off, the
vertical bar denotes the origin in $\Z$, and the arrow denotes the
position of the cursor. \label{fig:lampex1}}
\end{figure}

As described in \cite{deadlamp}, we consider an element in
 $L = \Z_2 \wr \Z$  geometrically as a particular configuration of a bi-infinite string
of light bulbs together with a ``lamplighter" or cursor.  Each
bulb has two states, on and off, and the cursor which indicates
the current bulb under consideration.  An element can be viewed as
a set of light bulbs which are illuminated, together with a
position of the cursor. Figure \ref{fig:lampex1} gives an example
of an element in $L$ represented in this way.

We consider a word in the wreath product generators $\{a,t\}$
representing an element $w \in L$.  We can view this word as a
sequence of instructions to create the configuration of
illuminated bulbs represented by the element $w$.  We begin with
the bi-infinite string of bulbs all  in the off state with the
cursor at the origin, then read the generators in the word one at
a time from left to right. The generator $a$ changes the state of
the bulb at the current cursor position,
 and the generators $t^{\pm1}$ move the cursor one position
right or left, respectively.  In this way, successive prefixes of
the word representing $w$ create a sequence of configurations of
light bulbs and cursor positions, ending with the one representing
$w$. We often view the cursor position as ``moving'' as dictated
by these prefixes, ending in the position corresponding to $w$. We
view the question of finding minimal length representatives for
$w$  as finding methods of constructing $w$ in this way as
efficiently as possible.

The wreath product generators $a$ and $t$ of $L$ encapsulate the
two instructions necessary to create an element $w \in L$. The
automata generators of $L$ combine these two basic operations of
moving the cursor and illuminating a bulb. The generator $ta$
moves the cursor one step to the right and changes the state of
the bulb at the new location of the cursor, and its inverse
$(ta)^{-1}= at^{-1}$ changes the state of the bulb at the current
location and then moves the cursor one step to the left.

The  location of the cursor in a word $w \in L$ is easily computed
as the exponent sum of $t$ if the word is described in terms of
the wreath product generating set.  For the automata generating
set,  the location of the cursor is the total exponent sum of $t$
and $ta$.  Since the relators have all have total $t$ and $ta$
exponent sum zero, the total exponent sums for different words
representing  the same element will be identical.

\subsection{Normal forms}
\label{sec:normalforms}

We would like to be able to calculate the word length of $w \in
L$ with respect to the automata generating set, in a manner
analogous to the
computation of word length with respect to the wreath product
generating set presented in \cite{deadlamp}.  This was
accomplished in \cite{deadlamp} through the use of normal forms
for elements of $L$.  Below, we use the same normal forms to aid
us in computing the word length of $w$ with respect to the
automata generating set.

Normal forms for elements of $L$ with respect to the wreath
product generators $a$ and $t$ are given in terms of the
conjugates $a_k=t^k a t^{-k}$, which move the cursor to the $k$-th
bulb, turn it on, and return the cursor to the origin.

We present two normal forms for an element $w \in L$ with
respect to the wreath product generators, the {\em right-first
normal form} given by
$$rf(w)=a_{i_1} a_{i_2} \ldots a_{i_m} a_{-j_1} a_{-j_2} \ldots a_{-j_l} t^{r}$$
and  the {\em left-first normal form} given by
$$lf(w)=a_{-j_1} a_{-j_2} \ldots a_{-j_l} a_{i_1} a_{i_2} \ldots a_{i_m}   t^{r}
$$
with $ i_m > \ldots i_2 > i_1 > 0 $ and $ j_l > \ldots j_2 > j_1
\geq  0 $. In the first $m$ prefixes of the right-first normal
form, only bulbs at positive indices are illuminated, in
increasing order.  In longer prefixes, the illuminated bulbs at
nonpositive indices are added as well.  Once our prefix includes
$a_{i_1} a_{i_2} \ldots a_{i_m} a_{-j_1} a_{-j_2} \ldots
a_{-j_l}$, all subsequent prefixes move the cursor closer its
position in $w$.

The left-first normal form has prefixes which first illuminate
bulbs at nonpositive integers. Note that these normal forms are
slightly different from the normal forms in \cite{deadlamp}, as
the bulb in position zero is grouped together with the bulbs at
negative integers rather than the bulbs at positive integers. The
reason for this grouping will become clear below.

One or possibly both of these normal forms will lead to a minimal
length representative for $w \in L$ with respect to the wreath
product generating set, depending upon the  location of the
cursor relative to the origin. That position is easy to detect
from the sign of the exponent sum of $t$, given as $r$ above. If
$r \geq 0$, then the left-first normal form will naturally lead to
a minimal length representative, and if $r \leq 0$, then the
right-first one will lead to a minimal length representative. If
the exponent sum of $t$ is zero, then both of these normal forms
yield minimal length representatives.  Words of this type, with
the cursor position at the origin, play an important role in the
geometry of the Cayley graph, as described in Section \ref{sec:cayley}

While these normal forms are defined in terms of conjugates of $a$
in wreath product
generating set, we use them to compute the word length of an
element with respect to the automata generating set as well.

We first note one difference between the two generating sets
discussed above, with consequences for computing the word length
of elements. In the wreath product generating set, travel and bulb
illumination are accomplished separately by the two generators.
That is, there must be at least one instance of the generator $a$
for each light bulb illuminated in the configuration of bulbs
representing a group element, and enough instances of the
generator $t$ to  move the cursor  to the appropriate bulbs.  In
the automata generating set, these two operations are combined in
the generator $ta$.  A single generator can either move the cursor
or move the cursor and illuminate a bulb in one step.

Second, we note the following about exponent sum in the prefixes
of a word representing $w \in L$.   As described above, each
prefix in a word representing $w$ determines a configuration of
illuminated bulbs, with the final configuration represented by
$w$.  Suppose that $w \in L$ has an illuminated bulb in position
$n$.  Then there must be a prefix of this word in which the cursor
is at position $n-1$.  We prove this in the following lemma.

\begin{lemma} \label{lem:leftlight}
Let $w \in L$ be represented by a word $\gamma$ in the
generators $\{t,ta\}$.  If $w$ has a bulb illuminated in position
$n$, there is some prefix of $\gamma$ for which the total exponent
sum of $ta$ and $t$ is $n-1$.
\end{lemma}

\begin{proof}
The only two generators which can illuminate a bulb are $ta$ and
$(ta)^{-1}=a t^{-1}$. Let $\eta$ be the shortest prefix of
$\gamma$ in which the bulb in position $n$ is illuminated. Recall
that the generator $t$ moves the cursor one unit to the right.  If
$\eta$ ends in the generator $ta$, then the prefix $\eta
(ta)^{-1}$ of $\gamma$ must have the cursor in position $n-1$;
that is, the exponent sum of all the generators in that prefix is $n-1$.

If $\eta$ ends in the generator $(ta)^{-1} = at^{-1}$, then the
cursor is at position $n-1$ in $\eta$, since after turning on the
bulb in position $n$ via the generator $a$, the generator $t^{-1}$
moves the cursor one position to the left. Thus, there is a prefix
of $\gamma$ with total exponent sum $n-1$.
\end{proof}

We again contrast the two generating sets.  With respect to the
wreath product generators, the cursor must be in position $-n$ to
illuminate a bulb in position $-n$, for $n \geq 0$.
With respect to the automata generating set, the cursor must move
to position $-n-1$ as part of illuminating a bulb in position
$-n$.  If this illumination is done via $ta$, then the cursor
begins in position $-n-1$ and moves to position $n$ to illuminate
the bulb there. If this illumination is done via $(ta)^{-1} =
at^{-1}$, then the cursor begins in position $-n$, illuminates the
bulb, and then moves to position $-n-1$.  Thus the automata
generating set requires the cursor to move farther from the origin
in order to illuminate the same bulb. With respect to illuminating
bulbs at positive indices, the need to visit bulb $n-1$ to illuminate
bulb $n$  does not result in any extra travel for the cursor, since
in that case the cursor will pass through position $n-1$ on its way
from the origin to position $n$.

We note that the length of $a$ with respect to the automata
generating set is two, as the cursor must be at position $-1$ at some
prefix of any word representing this element to ensure that the
bulb in position zero is illuminated. The minimal length
representatives of $a$ in this generating set are $t^{-1}(ta)$ and
$(ta)^{-1} t$.

To compute the length of an element with respect to the automata
generating set, the only data we require are the position of the
leftmost and rightmost illuminated bulbs  and the  cursor
position.  The configuration of the intermediate bulbs does not
affect the word length of an element.  In the prefixes of minimal
length representatives,  the cursor moves toward the extreme illuminated positions, then to its
position in the element.  At each position,
we choose $t$ or  $ta$ appropriately
to ensure that each bulb is left in the desired state. The  information about left-most and right-most
illuminated bulbs and cursor position can easily
be obtained from either normal form given above with respect to
the wreath product generating set, and thus we do not define new
normal forms for the automata generating set.

With these three pieces of information, it is easy to see how to
define minimal length representatives for elements of $L$ with
respect to the automata generators, in the spirit of the left- and
right-first normal forms.  If the cursor position in $w
\in L$ is to the right of the origin, we follow the example of
the left-first normal form and create a word whose initial
prefixes contain only the generators $t^{-1}$ and $(ta)^{-1}$, and
successively illuminate the correct bulbs to the left of the
origin.  Then we add enough instances of the generator $t$ to
create a prefix with the cursor at the origin.  We can then add a suffix to
this word containing only the generators $t$ and $ta$ to
illuminate the appropriate bulbs in positive positions and then finally
move the cursor to its position in $w$ with $t$ and $t^{-1}$ as
appropriate.
We will see below that this path gives a  minimal
length representative for $w$.  If the cursor position in $w$ is
at or to the left of the origin, we follow the example of the
right-first normal form and create a minimal length representative
in that way.  These representatives are not unique, as discussed below.

\subsection{Computing word length in $L$}
In \cite{deadlamp}, we used the normal forms given above to
compute the word length of $w \in L$ with respect to the wreath
product generating set as follows.
\begin{proposition}[\cite{deadlamp},Proposition 3.6]
Let $w = a_{i_1} a_{i_2} \ldots a_{i_m} a_{-j_1} a_{-j_2} \ldots
a_{-j_l} t^{r} \in L$, with \\ $0 < i_1 < i_2 \cdots < i_m $ and
$0 \leq j_1 < j_2 \cdots < j_l$. We define
$$D(w)=m+l+ min\{2 j_l+i_m + | r-i_m|,\  2 i_m+j_l+|r+j_l|\}.$$
The word length of $w$ with respect to the generating set
$\{a,t\}$ is exactly $D(w)$.
\end{proposition}

Thus, word length with respect to the wreath product generating set
is the number of illuminated bulbs  plus  the total necessary travel of
the cursor.

We now define a similar quantity, $D'(w)$,  also
computed from either normal form of $w$, which will give the word length of $w$
with respect to the automata generating set.

\begin{definition}
Let $w = a_{i_1} a_{i_2} \ldots a_{i_m} a_{-j_1} a_{-j_2} \ldots
a_{-j_l} t^{r} \in L$, with $0 < i_1 < i_2 \cdots < i_m $ and $0
\leq j_1 < j_2 \cdots < j_l$. If $l=0$, there are no bulbs
illuminated at or to the left of the origin and we set $D'(w)=
i_m+ | r -i_m|$. Otherwise, we set
 $$D'(w)= \min\{ 2 (j_l+1)+
 i_m+|r-i_m|,2 i_m+ j_l + 1 + |r+j_l+1|\}.$$
\end{definition}

We now prove that the quantity $D'(w)$, which is computed from the
normal form of $w$ in the wreath product generating set, yields
the word length of $w$ with respect to the automata generating
set.

\begin{proposition}
\label{prop:D} The word length of $w \in L$ with respect to
 the automata generating set
$\{t,ta\}$ is given by $D'(w)$.
\end{proposition}

Thus,  word length with respect to the automata generating set is merely the total travel required for the cursor, including the necessary travel to the bulb one to the left of the leftmost illuminated bulb.
Proposition \ref{prop:D} is proved via the following two lemmas.

\begin{lemma}
\label{lem:upper} The length of an element $w \in L$  with
 with respect to the generating set
$\{t,ta\}$ is at most $D'(w)$.
\end{lemma}

\begin{proof}
We first describe a word in $t$ and $ta$ representing the element
$w$ with the generators arranged in the order suggested by the
left-first normal form.  That is, we first have a sequence of
$j_l+1$ occurrences of the $t^{-1}$ or $a t^{-1}$ generators whose
prefixes move the cursor to the position $-j_l-1$, illuminating
the appropriate bulbs.  We add to this initial string of
generators $t^{j_l+1}$, so that this word now represents a
configuration of illuminated bulbs at non-positive positions with
the cursor at the origin.  We then add a suffix in the generators
$t$ and $ta$ of length $i_m$, so that the prefixes of this suffix
move the cursor to position $i_m$, illuminating the appropriate
bulbs at positive indices in increasing order. Finally, there are
$|r-i_m|$ occurrences of $t$ or $t^{-1}$ so that in the combined
word, the cursor is at the required position for the element $w$.
The length of this word is the first of the two terms over which
the minimum for $D'(w)$ is taken. Expressing $w$ as a word in the
manner suggested by the right-first normal form gives the second
of those two terms and thus we have the desired upper bound.

It is easy to see that if $l=0$,  there are no bulbs illuminated
in positions at or to the left of the origin.  In this case, we can avoid moving
the cursor to the left of the origin in the prefixes and we see that a similar argument to
above produces a word representing $w$ of length $i_m +
|r-i_m| = D'(w)$.
\end{proof}

\begin{lemma}
\label{lem:lower} The length of an element  $w \in L$  with
respect to the generating set $\{t,ta\}$ is at least $D'(w)$.
\end{lemma}

\begin{proof}
To compute a lower bound on the word length of $w \in L$, we
consider the minimum number of generators necessary to move the
cursor to the positions of the leftmost and rightmost illuminated
bulbs, and leave it in the proscribed position. As pointed out in
Lemma \ref{lem:leftlight}, to illuminate the bulb in position $n$,
the cursor must visit position $n-1$, which will result in additional
cursor travel in the prefixes of a representative in the case where there is an illuminated bulb
at or to the left of the origin.

We now compute this lower bound explicitly.  First, we  suppose that
$m>0$ and $l>0$, so that there is at least one illuminated bulb to the right of
the origin and at least one illuminated bulb at the origin or to the left of the origin.  We consider a word $\gamma$ in $\{t^{\pm
1},ta^{\pm 1}\}$ representing $w$ with specific prefixes.
We know that the bulbs in positions $i_m$ and $-j_l$ are
illuminated in $w$, and that the cursor position is $r \in \Z$.
Thus $\gamma$ must have prefixes with exponent sum, in order, $0,
i_m, -j_l-1, r$ or $0, -j_l-1, i_m, r$.  To accomplish this first
possible order, there must be at least $i_m + i_m + j_l +1 +
|r+j_l+1|$ occurrences of the generators, and to accomplish this
in the second order, there must be at least $j_l+1 + j_l+1 + i_m +
| r-i_m|$ occurrences of the generators, giving the desired
bounds.

If $m=0$, there are no bulbs illuminated to the right of the
origin, and the word $\gamma$ must have prefixes with exponent
sums $0, -j_l-1$, and $r$, giving a lower bound of $ j_l +1 +
|r+j_l+1|$ on the word length of $w$.

If $l=0$, there are no bulbs illuminated at the origin or to the
left of the origin, and the word $\gamma$ must have prefixes with
exponent sums $0, i_m$, and $r$, giving a lower bound of $ i_m + |
i_m-r|$ on the word length of $w$.
\end{proof}

Combining these lemma proves  Proposition \ref{prop:D} that the length of $w$ with respect to the
automata generating set $\{t,ta\}$ is exactly $D'(w)$.

As with the wreath product generating set, minimal length
representatives for group elements in $L$ are not generally
unique. If the total $t$ and $ta$ exponent sum in a word $w$ is
zero, the cursor is left at the origin in $w$ and there will be
minimal length representatives for $w$ arising from both the
left-first and right-first manners of construction described
above.  If there is any bulb which is ``visited twice'' by the
cursor during the construction of the element, there also will be
more than one minimal length representative for that element. That
is, if $\gamma$ is a minimal length representative for $w \in L$,
and $\gamma$ has two different prefixes $\gamma_1$ and $\gamma_2$
for which the  cursor position is $n$, then we can create another
minimal length representative for $w$.  If the bulb in position
$n$ is off in $w$, then we can either have it remain off during
the construction of $w$ or switch it on in $\gamma_1$ then later
off in $\gamma_2$. If the bulb in position $n$ is on in $w$, we
again have a choice as to which prefix illuminates  the bulb- it
can be switched on in $\gamma_1$ and remain on in $\gamma_2$ or be
left off in $\gamma_1$ and then switched on in $\gamma_2$.

For example, we consider the word $w=a_4 a_5 a_6 a_{-1} a_{-6}
t^{-2}$ pictured in Figure \ref{fig:lampex1}.  Since this word has
the cursor in position $-2$, to the left of the origin, only words
arising from the right-first normal form will be minimal and the
length of any minimal length representative is 24.  The bulb in
position $1$ will be visited twice during the construction of this
element, as will in fact all bulbs except those in positions $0$
and  $-1$.  Focusing on the bulb in position 1, we see that there
is a choice as to how many times the state of this bulb changes-
either twice or never,
 since it is left off in $w$.  So there is a minimal length representative
$t^3 (ta)^3 t^{-7} (ta)^{-1} t^{-4} (ta)^{-1} t^5$ of $w$ which
never turns on the bulb in position 1, and there is also the
minimal length representative $ (ta) t^2  (ta)^3 t^{-5} (ta)^{-1}
t^{-1} (ta)^{-1} t^{-4} (ta)^{-1} t^5$ of $w$ which switches the
bulb in position 1 on during the first visit and switches it off
during the second visit. We will have similar choices for bulbs in
positions 2, 3, 4, 5 and 6, as well as all bulbs in positions from
$-2$ to $-5$, resulting in many possible geodesic representatives
for $w$.

This analysis shows there are at least $2^u$ possible geodesic
representatives for group elements which have $u$ pairs of different
prefixes yielding the same cursor position and thus $u$ bulbs which are
visited twice by the cursor during the construction of the word.
  Furthermore, bulbs with the  cursor position at the origin will
have $2^u$ geodesics in each of the left-first and right-first directions,
giving $2^{u+1}$ total.

\section{Properties of the Cayley graph of $L$ as an automata group}
\label{sec:cayley}

\subsection{Dead-end elements}

We found in \cite{deadlamp} that $L$ contained elements which we
called {\em dead-end elements} with respect to the wreath product
generating set, meaning that a geodesic ray from the identity to
such an element could not be extended further.  We show that $L$
also contains dead-end elements with respect to the automata
generating set $\{t,ta\}.$

\begin{definition}
An element $w$ in a finitely generated group $G$ is a {\em dead-end
element with respect to a finite generating set $X$} for $G$
if $|w|=n$ and $|w x| \leq n$ for all generators $x$ in $X \cup
X^{-1}$, where $| \cdot |$ represents word length with respect to
the generating set $X$.
\end{definition}

These elements are called dead-end elements because a geodesic ray
in the Cayley graph $\Gamma(G,X)$ from the origin to a dead-end
element $w$ cannot be extended beyond $w$. Note that in groups
such as $L$ (with either generating set under consideration)
where all relators are of even length, if $w$ is a dead-end
element and $x \in X$, the word length of $wx$ will be necessarily
$n-1$.

There are different ``strengths" of dead-end behavior, measured by
the notion of {\em depth}.

\begin{definition}
An element $w$ in a finitely generated group $G$ is a {\em dead-end
element of depth $k$ with respect to a finite generating set $X$}
if $k$ is the largest integer with the following property: if
the word length of $w$ in $n$, then $|w x_1 x_2 \ldots x_l| \leq
n$ for $1 \leq l \leq k$ and all choices of generators $x_i \in X
\cup X^{-1}$.
\end{definition}

A dead-end element in the Cayley graph $\Gamma(G,X)$ is a point
from which is impossible to make immediate progress away from the
identity. The depth of a dead-end element $w$ reflects how long a
path from $w$ must be to reach an element further away from the
identity than $w$.

We show in \cite{ctcomb} that all dead-end elements in Thompson's
group $F$ with respect to the standard finite generating set
$\{x_0,x_1\}$ have depth two, and in \cite{deadlamp} that there
are dead-end elements of arbitrary depth in $L$ with respect to
the wreath product generating set.  The results of \cite{deadlamp}
extend to a larger class of wreath products as well.  We now show
that these results can be extended to the automata generating set
of $L$.

\begin{theorem}
\label{thm:deadends}
The lamplighter group $L$  contains dead-end
words of arbitrary depth  with respect to the generating set
$\{t,ta\}$.
\end{theorem}

\begin{figure}\includegraphics[width=3in]{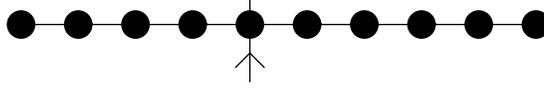}\\
\caption{The dead end word $d_5$, with all bulbs in positions -4 to 5
illuminated and the cursor at the origin. \label{fig:d5}}
\end{figure}

\begin{proof}
We define a family of elements which we show to be dead-end
elements with respect to the automata generating set.  Let $d_m$
denote a group element which has the bulbs at positions $m$ and
$-m+1$ illuminated, and the cursor at the origin. No bulbs may be
illuminated beyond positions $m$ and $-m+1$, and the state of the
intermediate bulbs between $m$ and $-m+1$ is irrelevant. In the example given
in Figure \ref{fig:d5}, we take all these intermediate bulbs  to be
illuminated yielding the element $d_m=a_0 a_1 a_2 \ldots a_m
a_{-1} a_{-2} \ldots a_{-m+1}$. These words all have length $4m$.
Since the cursor is at the origin, we can find minimal length
representatives for these words in the manner of either the left-first or right-first normal forms, as described in Section
\ref{sec:normalforms}.

We now check that these are dead-end elements.

\begin{itemize}
\item $w(t)$ and $w(ta)$ will have the cursor at position 1, and
their word length will now be $4m-1$. \item $w(t^{-1})$ and
$w(at^{-1})$ will have the cursor at position -1, and the total
length will again be $4m-1$.

\end{itemize}

To see that the depth of $d_m$ is at least $m$, we note the
following. An element with any possible configuration of bulbs
illuminated only between positions $-m+1$ and $m$, and with the
cursor remaining anywhere between $-m+1$ and $m$, will lie in the
$4m$ ball around the identity.  Thus, a path from $d_m$ to any
point in the $4m+1$ ball will have length at least $m+1$, since
the cursor must leave the range $[-m,m]$ in order to illuminate a
bulb with an index outside this range.
\end{proof}

We note that all elements of $L$ with the cursor at the origin and
illuminated bulbs to the right and at or to the left of the origin
are dead-end elements. The argument above shows that the depth is
at least the distance to the closer of the rightmost illuminated
bulb or one to the left of the leftmost illuminated bulb. In fact,
since the cursor motions to move to the rightmost or leftmost
illuminated bulbs all reduce word length, we see that the depth is
in fact twice the smaller of these two distances.  The specific
words $d_n$ used in the proof above to exhibit dead-end elements
of arbitrary depth have depth $2n$.

\subsection{Seesaw elements}

 An element $w$ in a finitely generated group $G$ is
a seesaw element with respect to a generating set $X$ if there are
only two possible suffixes that a geodesic representative for $w$
can have,  of the form $g^{\pm k}$ for a group generator $g$.
Phrased in terms of word length, we say that $w$ is a seesaw
element if there is a unique generator and its inverse which
reduce the word length of $w$, and for which $k$ subsequent
reductions in word length only occur through successive
applications of that generator. These elements may present
difficulty for the construction of canonical minimal length
representatives for group elements.   In \cite{seesaw} we show
that Thompson's group $F$ in the standard finite generating set
$\{x_0,x_1\}$ contains seesaw elements, and use them to show that
$F$ is not combable by geodesics. In \cite{deadlamp} we show that
the lamplighter group $L$ in the wreath product generating set as
well as a wide class of wreath products contain seesaw elements
with respect to at least one generating set.

We now precisely define seesaw elements.

\begin{definition}\label{def:seesaw}
A element $w$ in a finitely generated group $G$ with finite
generating set $X$ is a {\em seesaw element of swing $k>0$ with respect
to a generator $g$} if the following conditions hold.  Let $|w|$
represent the word length of $w$ with respect to the generating
set $X$.
\begin{enumerate}
\item  Right multiplication by both $g$ and $g^{-1}$ reduces the
word length of $w$; that is, $|wg^{\pm 1}| = |w| - 1$, and for all
$h \in X  \smallsetminus \left\{g^{\pm1}\right\}$, we have
$|wh^{\pm 1}| \geq |w|$.

\item Additionally, $|wg^l| = |wg^{l-1}| - 1$   for integral $l
\in [1,k]$, and $|wg^{m}h^{\pm 1}| \geq |wg^{m}|$ for all $h
\in X  \smallsetminus \left\{g\right\}$ and integral $m \in
[1,k-1]$.

\item Similarly,  $|wg^{-l}| = |wg^{-l+1}| - 1$  for integral $l
\in [1,k]$, and $|wg^{-m}h^{\pm 1}| \geq |wg^{-m}|$ for all $h
\in X  \smallsetminus \left\{g^{-1}\right\}$  for integral $m \in
[1, k-1]$.
\end{enumerate}
\end{definition}

These are called seesaw elements because they behave like a seesaw resting in
balance.  When in balance, there is a choice about which of two opposite directions to go down,
 but once that initial choice is made, there we can continue further downward only in that chosen
 direction, for the number of steps described as the swing.

 \begin{theorem}
\label{thm:noseesaw} The lamplighter group $L$  with respect to
the automata generating set does not contain any seesaw elements.
\end{theorem}

\begin{proof}
Suppose $w \in L$ is a seesaw element with respect to a generator
$g$.  Then both $g$ and $g^{-1}$ must reduce the word length of
$w$. If $g$ is $t$ or $t^{-1}$, then in order for both $t$ and
$t^{-1}$ to reduce word length, the cursor position in $w$ must be
at the origin with at least one bulb illuminated to the right of
the origin and at least one bulb illuminated at the origin or to
the left of the origin.  In this case, both $ta$ and $at^{-1}$
will also reduce word length, and $w$ cannot be a seesaw word as
the first condition of Definition \ref{def:seesaw} is not
satisfied.

Similarly, if $w$ is a seesaw element with respect to $ta$ or
$(ta)^{-1}$, in order for both $ta$ and $(ta)^{-1}$ to reduce word
length, the position of the cursor in $w$ must again be the
origin, with at least one bulb illuminated to the right of the
origin and at least one bulb illuminated at the origin or to the
left of the origin. In this case, both $t$ and $t^{-1}$ will also
reduce word length, and $w$ cannot be a seesaw word as the first
condition of Definition \ref{def:seesaw} is not satisfied.

Thus, for any $w \in L$ and $g \in \{t^{\pm 1},(ta)^{\pm 1}\}$, if
$w$ is a seesaw word then the position of the cursor in $w$ must
be the origin, and there must be bulbs illuminated at both
positive and non-positive indices.  However, in such a word
multiplication by all four generators will decrease the word
length.  Thus $w$ cannot be a seesaw word as the first condition
of Definition \ref{def:seesaw} is not satisfied.
\end{proof}

\subsection{Seesaw-like elements}

While $L$ with respect to the automata generating set does not
contain seesaw elements as it does with respect to the generating
set $\{a,t\}$, there are elements which exhibit similar behavior.
Seesaw words can be multiplied by a unique pair $g^{\pm 1}$ in
order to reduce their word length.  For the `seesaw-like' elements
we describe below, there are two distinct families of generators
analogous to the role of the generators $g$ and $g^{-1}$ for
seesaw words.

%

\begin{figure}\includegraphics[width=3in]{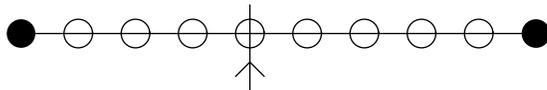}\\
\caption{The seesaw-like element $w_5= a_{-4} a_5$, which has the
cursor at the origin and length 20 with respect to the automata
generating set $\{t,ta\}$. \label{lampw5}}
\end{figure}

\begin{theorem}
\label{thm:seesaw} The lamplighter group $L$  with respect to the
automata generating set contains {\em seesaw-like} elements $w_k$,
which satisfy the following conditions with respect to the two
sets of generators $\{t,ta\}$ and $\{t^{-1},at^{-1}\}$.
\begin{enumerate}
\item   All generators reduce the length of $w_k$, that is $|w_k
g^{\pm 1}| = |w_k| -1 $ for all $g \in \{t, ta\}$.

\item Additionally, successive multiplication by the generators
$\{t, ta\}$ reduces word length for up to $k$ iterations.  That
is, $|w_k g_1 g_2 \ldots g_l| = |w_k | - l $   for integral $l \in
[1,k]$ and $g_i \in \{t, ta\} $. Furthermore,  any multiplication
by the other generators to these shortened words increases word
length, so we have  $|w_k g_1 g_2 \ldots g_l h| = |w_k | - l +1 $
for integral $l \in [1,k-1]$, $g_i \in \{t, ta\} $, and  $h \in
\{t^{-1}, a t^{-1}\} $.

\item Similarly, successive multiplication by the generators
$\{t^{-1}, a t^{-1}\}$ reduces word length for up to $k$
iterations and we have $|w_k h_1 h_2 \ldots h_l| = |w_k | - l $
for integral $l \in [1,k]$ and $h_i \in \{t^{-1}, a t^{-1}\} $.
Again, any multiplication by the other generators to these
shortened words increases word length, so we have  $|w_k h_1 h_2
\ldots h_l g| = |w_k | - l +1 $ for integral $l \in [1,k-1]$, $g
\in \{t, ta\} $, and  $h_i \in \{t^{-1}, a t^{-1}\} $.
\end{enumerate}
\end{theorem}


\begin{proof}
We consider the group elements $w_k= a_k a_{-k+1}$, which  have
length $4k$. The element $w_5$ is pictured in Figure \ref{lampw5}.
For a given $k$, this word has two illuminated bulbs in positions
$k$ and $-k+1$, and the cursor at the origin.

Since the cursor is at the origin in $w_k$, all generators reduce
word length and the element is a dead-end element.  Let $g \in
\{t^{\pm 1}, (ta)^{\pm 1}\}$.  Then in $w_kg$, the cursor is
either to the right or left of the origin in $\Z$.  It is clear
from the definition of $D'$ and Proposition \ref{prop:D} that once
the cursor is to the right of the origin, an application of $t$ or
$ta$ will decrease $D'$ and thus will  continue to reduce word
length until the position of the rightmost illuminated bulb is
reached. So for $l \in [1,k]$, the  group element $w_k g_1 g_2
\ldots g_l$ with $g_i \in \{t, ta\} $ has word length $|w_k|-l$.
Applications of $t^{-1}$ and $at^{-1}$ increase the word length of
these elements by moving the cursor to the left, except possibly
for the last one with $l=k$.  Similarly, a group element $w_k g_1
g_2 \ldots g_l$ with $g_i \in \{t^{-1}, (ta)^{-1}\} $ has word
length $|w_k|-l$, and computing $D'$ again shows that applications
of $t$ and $ta$ to those shortened words increase the word length
of these elements.
\end{proof}

Although these elements $w_k$ are not seesaw elements in the sense
of \cite{seesaw}, they do share an important property with seesaw
elements. In a true seesaw element, there is a choice between a
single generator and its inverse which both reduce word
length. After that initial choice, there are no further options
for reducing word length for the next $k$ steps; only repeated
applications of that same generator will reduce the word length to
$|w|-k$. Thus, this initial choice has long term consequences for
continued length reduction.

In the words described above, the choice for decreasing word
length is not between two generators, but two sets of generators:
$\{t,ta\}$ and $\{t^{-1},(ta)^{-1}\}$.  Each of these sets moves
the cursor position in a different direction, right and left,
respectively.  Once one of these sets of generators is chosen to decrease the
word length of a seesaw-like word $w_k$, only generators from that
{\em set} of generators can further decrease the word length for
$k$  iterations.   As in true seesaw words,  these two different directions
diverge as quickly as possible,  though
both directions reduce word length.  For example, for $s \leq k$,  the distance between the shortened
elements $w_k t^s$ and $w_k t^{-s}$ is $2s$.

\bibliographystyle{plain}

\end{document}